\documentclass[a4paper,11pt, leqno]{article}
\RequirePackage[intlimits]{amsmath}
\RequirePackage{amssymb}
\RequirePackage{amsthm}
\RequirePackage{latexsym}
\RequirePackage[mathscr]{eucal}
\RequirePackage{ifthen}
\RequirePackage{verbatim}

\newcounter{segno}[section]

\newskip\segskipamount
\newskip\procskipamount
\newskip{\interskipamount}
\newcommand{\segskip}{\vskip\segskipamount}
\newcommand{\procskip}{\vskip\procskipamount}
\newcommand{\interskip}{\vskip\interskipamount}

\renewcommand{\thesegno}{\thesection.\arabic{segno}}

\newenvironment{segment}{%
  \segskip
  \refstepcounter{segno}%
  \noindent\textbf{(\thesegno)\ }}%
  {\segskip}

\newenvironment{segproc}[1]{%
  \textbf{#1.\ }%
  \em}%
  {}

\newenvironment{proclamation}[1]{%
  \par\procskip
  \noindent\textbf{#1.\ }%
  \em}%
  {}

  {}%

\numberwithin{equation}{segno}

\makeatletter
\newenvironment{assertionlist}{%
  \renewcommand{\theenumi}{\arabic{enumi}}%
  \renewcommand{\theenumii}{\roman{enumii}}%
  \renewcommand{\p@enumii}{\theenumi.}%
  \renewcommand{\p@enumiii}{\theenumi.\theenumii.}%
  \begin{enumerate}}
  {\end{enumerate}}
\makeatother

\makeatletter
\newenvironment{definitionlist}{%
  \renewcommand{\theenumi}{\alph{enumi}}%
  \renewcommand{\theenumii}{\roman{enumii}}%
  \renewcommand{\p@enumii}{\theenumi.}%
  \begin{enumerate}}
  {\end{enumerate}}
\makeatother

\makeatletter
\newenvironment{bulletlist}{%
  \renewcommand{\theenumi}{$\bullet$}%
  \begin{enumerate}}
  {\end{enumerate}}
\makeatother

\RequirePackage[intlimits]{amsmath}
\RequirePackage{amssymb}
\RequirePackage{latexsym}
\RequirePackage[mathscr]{eucal}
\RequirePackage[scriptlabels]{diagrams}
\RequirePackage{ifthen}
\RequirePackage{verbatim}

\newcommand{\FF}{\mathbb{F}}
\newcommand{\GG}{\mathbb{G}}

\newcommand{\QQ}{\mathbb{Q}}

\newcommand{\ZZ}{\mathbb{Z}}

\newcommand{\kbar}{\bar{k}}

\newcommand{\Gbar}{\bar{G}}

\newcommand{\Zbar}{\bar{Z}}

\newcommand{\zgbar}{\bar{\zeta}}

\newcommand{\Ihat}{\hat{I}}

\newcommand{\What}{\hat{W}}

\newcommand{\gtilde}{\tilde{g}}

\newcommand{\Ltilde}{\tilde{L}}
\newcommand{\Mtilde}{\tilde{M}}

\newcommand{\Ptilde}{\tilde{P}}

\newcommand{\Rtilde}{\tilde{R}}

\newcommand{\Ytilde}{\tilde{Y}}
\newcommand{\Ztilde}{\tilde{Z}}

\newcommand{\mfr}{{\mathfrak m}}

\newcommand{\Ascr}{{\mathscr A}}

\newcommand{\Gscr}{{\mathscr G}}
\newcommand{\Hscr}{{\mathscr H}}

\newcommand{\Oscr}{{\mathscr O}}
\newcommand{\Pscr}{{\mathscr P}}

\newcommand{\Uscr}{{\mathscr U}}

\newcommand{\Wscr}{{\mathscr W}}

\newcommand{\Ascrtilde}{\tilde\Ascr}

\newcommand{\eps}{\varepsilon}

\newcommand{\agtilde}{\tilde\alpha}

\newcommand{\zgtilde}{\tilde\zeta}

\DeclareMathOperator{\Bor}{Bor}
\DeclareMathOperator{\Cent}{Cent}

\DeclareMathOperator{\codim}{codim}

\DeclareMathOperator{\DR}{DR}

\DeclareMathOperator{\gr}{gr}

\DeclareMathOperator{\Isom}{Isom}

\DeclareMathOperator{\Par}{Par}

\DeclareMathOperator{\relpos}{relpos}

\DeclareMathOperator{\rk}{rk}
\DeclareMathOperator{\Sp}{Sp}
\DeclareMathOperator{\Spec}{Spec}

\newcommand{\dlbrack}{\mathopen{[\![}}
\newcommand{\drbrack}{\mathclose{]\!]}}
\newcommand{\dual}{{{}^{\vee}}}
\newcommand{\Isomline}{\underline{\Isom}}

\newcommand{\lrangle}{{\langle\ ,\ \rangle}}

\newcommand{\restricted}[1]{\vert_{#1}}
\newcommand{\set}[2]{\{\,#1 \mid \text{#2}\,\}}

\newarrow{To}---->
\newarrow{Epi}----{>>}
\newarrow{Mono}{boldhook}--->
\newarrow{Bij}<--->
\newarrow{Sends}|--->

\newcommand{\ariso}{\overset{\sim}{\rightarrow}}

\newcommand{\air}{\hookrightarrow}

\DeclareMathOperator{\SP}{SP}
\DeclareMathOperator{\GSp}{GSp}

\begin{document}

\title{Flatness of the mod $p$ period morphism for the moduli space
  of principally polarized abelian varieties}

\author{Torsten Wedhorn}


\maketitle

\noindent{\scshape Abstract.}
In this note it is shown that associating to a principally polarized
abelian variety its de Rham cohomology defines a faithfully flat
morphism of the moduli space of principally polarized abelian
varieties in positive characteristics to the moduli space of
symplectic $F$-Zips as defined in \cite{MW}.

\bigskip

\section*{Introduction}

Let $\Ascr_g$ be the moduli space of $g$-dimensional principally
polarized abelian varieties in characteristic $p > 0$. If $S$ is any scheme of
characteristic $p$, an $S$-valued point consists of an abelian scheme
$f\colon A\to S$
of dimension $g$ and a principal polarization $\lambda$ of $A$. We can
associate to this pair the following ``linear algebra'' datum over
$S$:

We set $M = H^1_{\DR}(A/S)$ which is a locally free
$\Oscr_S$-module of rank $2g$. The polarization $\lambda$ induces a
perfect alternating bilinear form on $M$.

There are two canonical spectral
sequences converging against the\break de~Rham cohomology, namely the Hodge
spectral sequence and the conjugate spectral sequence. Both of them
induce a filtration on $M$. Moreover, they both degenerate and the
filtrations are of the form
\begin{gather*}
0 \rTo C := f_*\Omega^1_{A/S} \rTo M \rTo R^1f_*\Omega^0_{A/S} \rTo
 0,\\
0 \rTo D := R^1f_*(\Hscr^0(\Omega^{\bullet}_{A/S})) \rTo M \rTo
f_*(\Hscr^1(\Omega^{\bullet}_{A/S})) \rTo 0.
\end{gather*}
Finally, the Cartier isomorphism induces isomorphisms
\[
\varphi_0\colon (M/C)^{(p)} \rTo^{\sim} D, \qquad \varphi_1\colon C^{(p)}
\rTo^{\sim} M/D.
\]

This tuple is a ``symplectic $F$-zip of type $(g,g)$'' in the terminology of
\cite{MW}. In other words, we have constructed a morphism $\zgbar$ from
$\Ascr_g$ into the moduli space $\Zbar_g$ of symplectic $F$-zips of
type $(g,g)$. We show:

\begin{proclamation}{Theorem}
The morphism $\zgbar$ is faithfully flat.
\end{proclamation}

\smallskip

In \cite{MW} it was proved that $\Zbar_g$ is the quotient (in the
sense of algebraic stacks) of a smooth variety $Z_J$ by the action of
the symplectic group $G = Sp_{2g}$ and that the $G$-orbits on $Z_J$ are in
natural bijection with elements in $W_J\backslash W$ where $W$ is the
Weyl group of $G$ and $W_J$ is a certain Levi subgroup of
$W$. Moreover, a formula for the codimension of each stratum was given.

In this particular case, $W_J\backslash W$ can be identified with
$\{0,1\}^g$ and the codimension of the $G$-orbit corresponding to
$u = (\varepsilon_i) \in \{0,1\}^g$ is equal to
\[
g(g+1)/2 - \sum_{i=1}^g \epsilon_i(g+1-i).
\]

Every $G$-orbit $Z^u_J$ in $Z_J$ defines a locally closed substack
$\Zbar^u_g$ of $\Zbar_g$. The inverse images of these substacks in
$\Ascr_g$ are the Ekedahl-Oort strata $\Ascr^u_g$ of $\Ascr_g$ as
defined in \cite{Oo}. As $\dim(\Ascr_g) = g(g+1)/2$, the flatness of
$\zgbar$ implies:

\begin{proclamation}{Corollary}
For $u = (\varepsilon_i) \in \{0,1\}^g$ let $\Ascr^u_g$ be the
corresponding Ekedahl-Oort stratum. Then $\Ascr^u_0$ is
equi-dimensional and we have
\[
\dim(\Ascr^u_g) = \sum_{i=1}^g \epsilon_i(g+1-i).
\]
\end{proclamation}

\smallskip

This has also been shown by Oort in \cite{Oo} by different methods
(the elementary sequence $\varphi$ in the terminology of loc.~cit.\ is
given by $\varphi(i) = \eps_1 + \dots + \eps_i$).

\medskip

We will now give a short overview over the structure of the article. In
the first section we fix some notations from the theory of reductive
groups. In section 2 the moduli space of principally polarized abelian
varieties is defined. Section 3 contains the development of the main technical
tool for the proof of the flatness. Here we define split Dieudonn\'e
displays. On one hand they are related to Dieudonn\'e displays defined
by Zink \cite{Zi} and we will use Zink's theory for our proof. On the
other hand we can express the datum of a split Dieudonn\'e display in
a purely group theoretical way (see \eqref{displaytogroup}). In the
last section we define the
morphism of (a cover of) our moduli spaces to (a cover of) the moduli
space of $F$-zips with additional structures and then show that this
morphism is flat and surjective.

\smallskip

\noindent{\itshape Notation}. Let $G$ be a group, $X \subset G$ a subset and
$g \in G$. Then we set ${}^gX = gXg^{-1}$. 


\section{Preliminaries on reductive groups}

\begin{segment}\label{defWeyl}
Let $k$ be a field and let $G$ be a connected reductive group over
$k$. We fix an algebraic closure $\kbar$ of $k$. We denote by $\Bor_G$
the scheme of Borel groups of $G$.

In $\Bor_G \times \Bor_G$ we define the subscheme $\SP$ whose
$S$-valued points consists of those pairs $(B,B')$ of Borel subgroups
of $G_S$ such that fppf-locally there exists a maximal torus $T$ of
$G$ which is contained in $B$ and $B'$. In this case we say that $B$
and $B'$ are in \emph{standard position}.

The group $G$ acts on $\SP$ by simultaneous conjugation and the
fppf-quotient of this action is representable by a finite \'etale
$k$-scheme $\Wscr_G$. We set $W := W_G := \Wscr(\kbar)$ and call it
the \emph{Weyl group} of $G$. For any $w
\in W$ we denote the corresponding $G(\kbar)$-orbit of $\SP(\kbar)$ by $\SP^w$.

The set $W$ can be endowed with a group structure: For
$w,w' \in W$ choose $(B_1,B_2) \in \SP^w$ and $(B'_1,B'_2) \in
\SP^{w'}$ such that $B_2 = B'_1$ and such that the three Borel
subgroups $B_1$, $B_2 = B'_1$ and $B'_2$ contain a common maximal
torus (this always can be done). Then we define
the product $ww'$ as the $G$-orbit of $(B_1,B'_2)$.

For any $w \in W$ we define the \emph{length} of $w$ as
\[
\ell(w) = \dim(\SP^w) - \dim(\Bor).
\]
We call the set $I := I_G := \set{w \in W}{$\ell(w) = 1$}$ the set of
\emph{simple elements} in $W$.

It is well known that $(W,I)$ is a Coxeter group and that $\ell$ is
the usual length function with respect to the Coxeter base $I$ (see
\cite{Lu}~7.3).

The canonical map $\SP(\kbar) \to W$ is denoted by $\relpos$ and
$\relpos(B,B')$ is called the \emph{relative position} of $B$ and $B'$.
\end{segment}

\begin{segment}
Let $\varphi\colon G \to G'$ be a homomorphism of reductive groups
such that the induced homomorphism $\varphi^{\rm ad}\colon G^{\rm ad}
\to G^{\prime{\rm ad}}$ of adjoint groups is an isomorphism. Then
$\varphi$ induces an isomorphism $(W_G,I_G) \ariso (W_{G'},I_{G'})$ of
Coxeter groups.

Indeed, it suffices to show that the canonical homomorphism
$G \to G^{\rm ad}$ induces an isomorphism
of Weyl groups. But for every $k$-scheme $S$ every Borel subgroup of
$G_S$ contains the center of $G_S$ and therefore $\varphi$ induces an
isomorphism $\Bor_G \ariso \Bor_{G^{\rm ad}}$. As the center of $G$
acts trivially on $\Bor_G$, the claim follows.
\end{segment}

\begin{segment}\label{parabolics}
We denote by $\Par = \Par_G$ the scheme of parabolic subgroups of
$G$. The algebraic group $G$ acts on $\Par$ via conjugation.
The $G_{\kbar}$-orbits of $\Par_{\kbar}$ are the connected components
of $\Par_{\kbar}$.

For $i \in I$ we say that $P \in \Par(\kbar)$ is of type $\{i\}$
if for all pairs $(B,B')$ of Borel subgroups $B \not= B'$ which are
contained in $P$ we have $\relpos(B,B') = i$.

If $P$ is any parabolic subgroup of $G_{\kbar}$, we define the type $J$
as the subset of $I$ which consists of those $i \in I$ such that $P$
contains a parabolic subgroup of type $\{i\}$. This sets up a
bijection between $G(\kbar)$-conjugacy classes in
$\Par(\kbar)$ and the set of subsets of $I$.

For every such subset $J$ we denote by $\Par_J$ the variety of
parabolics of type $J$. This variety is defined over the field
extension of $k$ in $\kbar$ over which $J \subset W = \Wscr(\kbar)$
is defined. We denote by $\Pscr_J$ the universal parabolic of type $J$
over $\Par_J$ and by $\Uscr_J$ its unipotent radical.
\end{segment}

\begin{segment}\label{defLeviWeyl}
For any subset $J$ we denote by $W_J$ the subgroup of $W$ generated by
$J$. Alternatively, $W_J$ can be defined as the set of $\relpos(B,B')$
where $(B,B')$ runs through all pairs of Borel subgroups which are
contained in a common parabolic subgroup of type $J$.

For $w \in W$ and $J,K \subset I$ the double coset $W_JwW_K$ contains
a unique element of minimal length and we set
\[
{}^JW^K = \set{w \in W}{$w$ is of minimal length in $W_JwW_K$}.
\]
We define ${}^JW = {}^JW^\emptyset$ and $W^K = {}^\emptyset W^K$.

For any two parabolic subgroups $P \in \Par_J$ and $Q \in \Par_K$ we
consider the set
\[
\set{\relpos(B,C)}{$B \subset P$, $C \subset Q$ Borel subgroups}
\subset W.
\]
Clearly $W_J$ acts from the left and $W_K$ acts from the right on this
set, and it consists only of a single orbit under this $(W_J, W_K)$-
action. Therefore, it contains a unique element of minimal length which
we call the relative position of $P$ and $Q$ and which we denote by
$\relpos(P,Q) \in {}^JW^K$.
\end{segment}

\begin{segment}\label{symplecticexample}
As an example we consider a symplectic space $(V,\lrangle)$ of
dimension $2g$ over a field $k$ and denote by $G = GSp(V,\lrangle)$
the group of symplectic similitudes of $(V,\lrangle)$. Moreover, we
denote by $G' = Sp(V,\lrangle) \subset G$ the group of symplectic
isomorphisms.

As $G/G' \cong \GG_m$ is abelian, we have an isomorphism $\Par_G \ariso
\Par_{G'}$ given on points by $P \mapsto P \cap G'$. This identifies
in particular $\Bor_G$ with $\Bor_{G'}$. As $G$ is the product of $G'$
and $\Cent(G)$, this isomorphism also induces an isomorphism
$\Wscr_{G'} \ariso \Wscr_{G}$. In the sequel we will only consider $G'$
and set $\Bor = \Bor_{G'}$, $\Wscr = \Wscr_{G'}$ and so on.

Let $S$ be any $k$-scheme. Then $V_S = V \otimes_k \Oscr_S$ is a free
$\Oscr_S$-module of rank $2g$, and the base change of $\lrangle$ is a
perfect alternating form on $V_S$.

Let $M$ be any locally free $\Oscr_S$-module locally of finite
rank. An $\Oscr_S$-submodule $F \subset M$ is called locally direct
summand if one of the following equivalent conditions holds:
\begin{assertionlist}
\item For every open affine subset $U$ of $S$ there exists a
complement $F'$ of $F\restricted{U}$ in $M\restricted{U}$ (i.e.~an
$\Oscr_U$-submodule $F'$ of $M\restricted{U}$ such that
$F\restricted{U} \oplus F' = M\restricted{U}$).

\item Locally for the fpqc-topology $F$ admits a complement in $M$.

\item The quotient $M/F$ is a locally free $\Oscr_S$-module.

\item There exists an open affine covering $(U_i = \Spec(A_i))_{i\in
    I}$ of $S$ such that for all $i \in I$ and every maximal ideal
    $\mfr$ of $A_i$ the base change 
\[
\iota \otimes \kappa(\mfr)\colon F_{\kappa(\mfr)} \to M_{\kappa(\mfr)}
\]
is injective where $\iota\colon F \to M$ is the inclusion
homomorphism.
\end{assertionlist}

If $F$ is locally a direct summand of $V_S$, its orthogonal dual
$F\dual \subset V^*_S$ is locally a
direct complement and therefore its orthogonal $F^{\bot} \subset V_S$
is locally a direct summand. Here we denote by $(\ )^*$ the
$\Oscr_S$-linear dual.

A \emph{flag} of $V_S$ is a sequence
\[
0 = F_0 \subsetneq F_1 \subsetneq \dots \subsetneq F_{r-1} \subsetneq F_r = V_S
\]
of $\Oscr_S$-submodules of $V_S$ which are all locally direct summands
of $V_S$. This implies that $F_i$ is also a locally direct summand of
$F_{i+1}$. Such a flag is called \emph{complete} if $F_{i+1}/F_i$ is
locally free of rank 1 and it is called \emph{symplectic} if for all $i$
there exists a $j$ such that $F_i^{\bot} = F_j$. Any symplectic flag
is already uniquely determined by those members which are totally isotropic.

We denote by $\Par'$ the $k$-scheme whose $S$-valued points are the
symplectic flags and by $\Bor'$ the subscheme whose $S$-valued points
are those which consist of complete symplectic flags. These are projective
and smooth $k$-schemes. By associating to every symplectic flag over $S$ its
stabilizer in $G'_S$ we obtain an isomorphism $\Par' \ariso \Par$ which
induces an isomorphism $\Bor' \ariso \Bor$.

We use this ismorphism to identify symplectic flags of $V_S$
and parabolic subgroups of $G'_S$. Two complete symplectic flags
$(F_i)$ and $(F'_i)$ are in standard position if and only if $F_i + F'_j$ is
locally a direct summand of $V_S$ for all $i,j = 1,\dots,2g$. By
passage to the orthogonal it suffices to check this condition for
those $F_i$ and $F'_j$ which are totally isotropic, i.e.~for $i,j =
1,\dots g$.

Two pairs of flags $((F_i),(G_i))$ and $((F'_i),(G'_i))$ of $V_S$ are in the
same $G'(S)$-orbit if the rank of $F_i + G_j$ equals the rank of $F'_i
+ G'_j$ for all $i,j = 1,\dots,g$. From this it is easily seen that
every $G(\kbar)$-orbit of a pair of complete symplectic flags in
$V_{\kbar}$ contains a pair of complete symplectic flags in $V$. Therefore,
the finite \'etale scheme $\Wscr$ is the constant scheme associated to
the set $W$.

We describe $W$: Let $((F_i),(G_j))$ be a pair of complete symplectic
flags in $V$. For all $i$ there is a unique $j = \pi(i)$ such that
$\gr^G_j\gr^F_i \not= (0)$. The rule $i \mapsto \pi(i)$ defines a
permutation of the set $\{1,\dots,2g\}$, i.e.~an element $\pi \in
S_{2g}$.  Clearly, two pairs of flags are in the same $G$-orbit if and
only if the associated permutations are equal. Moreover, the condition
for the flags to be symplectic implies that the permutation
satisfies the condition
\begin{alignat}{2}\label{symplecticperm}
\pi(i) + \pi(2g+1-i) &= 2g+1 & & \qquad\text{for all $i = 1,\dots,g$.}
\end{alignat}
 Hence we can identify $W$ with the set of $\pi \in S_{2g}$ satisfying
 \eqref{symplecticperm}. One can check that this defines a group
 isomorphism.

The simple elements in $W$ correspond to those $G$-orbits of pairs
of symplectic complete flags $((F_i),(G_j))$ such that there exists an
$i_0 \in \{1,\dots,g\}$ with $F_{i_0} \not= G_{i_0}$ and $F_i = G_i$ for all $i
\in \{1,\dots,g\}$ with $i \not= i_0$. As we have $F_j^{\bot} =
F_{2g-j}$ for all $j$, this implies $F_{2g-i_0} \not= G_{2g-i_0}$ and
$F_i = G_i$ for all $i \in \{1,\dots,2g\} \setminus
\{i_0,2g-i_0\}$. Therefore, the set of simple elements $I$ consists of
$\{s_1,\dots,s_g\}$ with
\[
s_i =
\begin{cases} \tau_i\tau_{2g-i}, &\text{for $i = 1,\dots,g-1$;}\\
\tau_g, &\text{for $i = g$.}
\end{cases}
\]
where $\tau_j \in S_{2g}$ denotes the transposition of $j$ and $j+1$.

Let $J = \{s_{i_1},\dots,s_{i_r}\}$ be a subset of $I$. If $P$ is a
parabolic of $G$ corresponding to a symplectic flag $(F_j)$, then $P$
is of type $J$ if and only if for all $\rho = 1,\dots,r$ and for all
$j$ the rank of $F_j$ is not equal to $i_{\rho}$.
\end{segment}

\begin{segment}\label{Siegelparabolic}
Consider the special case $J = \{s_1,\dots,s_{g-1}\}$. Then $W_J$
consists of those permutation $\pi \in W$ such that
$\pi(\{1,\dots,g\}) = \{1,\dots,g\}$. The map
\[
W_J \to S_g. \qquad \pi \mapsto \pi\restricted{\{1,\dots,g\}}
\]
is a group isomorphism. A symplectic flag is of type $J$
if and only if it is of the form $(0) \subset F \subset V_S$ where $F$
is locally a direct summand with $F^{\bot} = F$.

The set ${}^JW$ consists in this case of those elements $\pi \in W$
such that
\[
\pi^{-1}(1) < \pi^{-1}(2) < \dots < \pi^{-1}(g).
\]
Of course, this implies $\pi^{-1}(g+1) < \dots < \pi^{-1}(2g)$. If
$\Sigma = \{j_1 < \dots < j_g\} \subset \{1,\dots,2g\}$ is a subset of
$g$ elements such that either $i \in \Sigma$ or $2g+1-i \in \Sigma$
for all $i = 1,\dots,g$, we get a corresponding element $\pi_{\Sigma}
\in {}^JW$ by setting $\pi^{-1}(i) = j_i$.

The sets of these $\Sigma$'s is in bijection with $\{0,1\}^g$ by
associating to $\Sigma$ the tuple $(\epsilon_1,\dots,\epsilon_g)$ with
\[
\epsilon_i =
\begin{cases}
0, &\text{if $i \in \Sigma$;}\\
1, &\text{otherwise.}
\end{cases}
\]
The length of such an element $(\epsilon_1,\dots,\epsilon_g)$ is equal
to
\[
\sum_{i=1}^g \epsilon_i(g+1-i).
\]
\end{segment}


\section{The moduli space}

\begin{segment}\label{PELDatum}
We fix a symplectic vector space $(V,\lrangle)$ over $\QQ$ and denote
by $2g$ its dimension. We assume that $g > 0$. Moreover, we fix a
prime $p > 0$. Let $\Lambda$ be a $\ZZ_p$-lattice in
$V_{\QQ_p}$ such that the restriction of $\lrangle_{\QQ_p}$ to
$\Lambda$ is perfect.

Let $G = \GSp(V,\lrangle)$ be the group of symplectic similitudes of
$(V,\lrangle)$ and denote by $G' \subset G$
the subgroup of symplectic isomorphisms. We consider $G$ and $G'$ as
reductive groups over $\QQ$.

We set $\Gscr = \GSp(\Lambda,\lrangle)$, $\Gscr' =
\Sp(\Lambda,\lrangle)$. These are reductive group schemes over $\ZZ_p$
whose generic fibres is equal to $G_{\QQ_p}$ and $G'_{\QQ_p}$,
respectively. Their special fibres are denoted by $\Gbar$ and
$\Gbar'$, respectively.
\end{segment}

\begin{segment}
We denote by $\Ascr = \Ascr_g$ the $\ZZ$-groupoid whose fibres
over a scheme $S$ consists of the category of tuples
$(A,\lambda)$ where
\begin{bulletlist}
\item
$A$ is an abelian scheme over $S$ of dimension $g$;
\item
$\lambda$ is a principal polarization of~$A$.
\end{bulletlist}

Then $\Ascr$ is a smooth algebraic Deligne-Mumford stack over $\ZZ$ of
relative dimension $g(g+1)/2$.
We denote by $\Ascr_0 = \Ascr_{g,0}$ the reduction
$\Ascr \otimes_{\ZZ} \FF_p$ at~$p$.
\end{segment}

\begin{segment}\label{associateFzip}
Let $S$ be a $\FF_p$-scheme and let
$(A,\lambda)$ be an $S$-valued point of~$\Ascr_0$.
We set $M = H^1_{\DR}(A/S)$. Then $M$ is a locally free
$\Oscr_S$-module of rank $2g$ endowed with a perfect
alternating form $\lrangle$ induced by $\lambda$. Let $\beta\colon M
\ariso M^* = \Hscr om_{\Oscr_S}(M,\Oscr_S)$ be the isomorphism associated
to $\lrangle$.

Moreover, $M$ has two locally direct summands, namely $C :=
f_*\Omega^1_{A/S}$ given by the Hodge spectral sequence and $D :=
R^1f_*(\Hscr^0(\Omega^{\bullet}_{A/S}))$ given by the conjugate
spectral sequence. Both $C$ and~$D$ are locally direct summands of rank
$g$. We have $D = D^{\bot}$ and~$C = C^{\bot}$, i.e.~$\beta$ induces
isomorphisms $M/C \ariso C^*$ and $C \ariso (M/C)^*$ (similarly for $D$).

Finally, the Cartier isomorphism induces isomorphisms
\[
\varphi_0\colon (M/C)^{(p)} \rTo^{\sim} D, \qquad \varphi_1\colon C^{(p)}
\rTo^{\sim} M/D
\]
such that the diagram
\[
\begin{diagram}
(M/C)^{(p)} & \rTo^{\varphi_0} & D \\
\dTo_{\beta^{(p)}} & & \dTo^{\beta} \\
(C^*)^{(p)} & \lTo^{\varphi_1^*} & (M/D)^* \\
\end{diagram}
\]
commutes.
\end{segment}

\begin{segment}\label{WeylShimura}
We denote by $(W,I)$ the Weyl group of $G$ together with its set of
simple reflections \eqref{defWeyl}. Let $J \subset I$ be the subset of
simple reflections corresponding to the conjugacy class of those
parabolic subgroups which are the stabilizers of symplectic flags of
the form $(0) \subsetneq F \subsetneq V_S$ \eqref{Siegelparabolic}.
As $G'$ and $G$ have the same adjoint group, the
Weyl group of $G'$ is canonically identified with the Weyl group
of $G$.

Via the reductive group scheme $\Gscr'$ we have a
canonical isomorphism of $(W,I)$ with the Weyl group of $\Gscr'
\otimes \FF_p = \Gbar'$.
\end{segment}

\begin{segment}\label{localisom}
\begin{segproc}{Lemma}
Let $T$ be any scheme and let $M_1$ and $M_2$ be two
locally free $\Oscr_T$-modules of the same rank with a symplectic form
$\lrangle_i$. Then Zariski locally on $T$, the symplectic modules
$(M_1,\lrangle_1)$ and $(M_2,\lrangle_2)$ are isomorphic. Moreover, the
scheme of symplectic isomorphisms
\[
\Isomline_{\Oscr_T}((M_1,\lrangle_1),(M_2,\lrangle_2))
\]
is smooth.
\end{segproc}

\begin{proof}
This is clear.
\end{proof}
\end{segment}

\begin{segment}\label{covermoduli}
We define two smooth coverings $\Ascr^{\#}_0$ and
$\Ascrtilde_0$ of~$\Ascr_0$ as follows:

For every
$\FF_p$-scheme~$S$ the $S$-valued points of $\Ascr^{\#}_0$
are given by tuples $(A,\lambda,\alpha)$ where
$(A,\lambda) \in \Ascr_0(S)$ and where $\alpha$ is
an $\Oscr_S$-linear symplectic isomorphism $H^1_{\DR}(A/S) \ariso
\Lambda \otimes_{\ZZ_p} \Oscr_S$.

Therefore, $\Ascr^{\#}_0$ is a torsor
for the Zariski topology over~$\Ascr_0$ under the smooth group
scheme~$\Gbar'$.

The $S$-valued points of $\Ascrtilde_0$ are given by tuples
$(A,\lambda,\alpha, C',D')$ with
$(A,\lambda,\alpha) \in \Ascr^{\#}_0$ and where
$C'$ and~$D'$ are totally isotropic
complements of $\alpha(C)$ and of~$\alpha(D)$ in $\Lambda \otimes_{\ZZ_p}
\Oscr_S$, respectively.
\end{segment}


\section{Dieudonn\'e displays with additional structures}

\begin{segment}\label{defWhat}
In this section we will always denote by $R$ a complete local
Noetherian ring $R$ with perfect residue field of characteristic
$p$. If $p = 2$, we also assume that $pR = 0$.

We will endow Dieudonn\'e displays in the sense of
Zink \cite{Zi} with additional structures. We use freely the
terminology of loc.~cit.. In particular, we have the ring $\What(R)$
which is endowed with Frobenius $\sigma$ and
Verschiebung $\tau$ (which are denoted by ${}^F$ and ${}^V$,
respectively, in loc.~cit.). The kernel of the canonical homomorphism
$\What(R) \to R$ is denoted by $\Ihat_R$. Note that we have $\What(k)
= W(k)$ if $k$ is a perfect field. There is a unique structure of a
$\ZZ_p$-algebra on $\What(R)$.

For every $\What(R)$-module $M$ we set
$M^{\sigma} =  \What(R) \otimes_{\sigma,\What(R)} M$.
If $M$ is of the form $M = \Lambda \otimes_{\ZZ_p} \What(R)$ for some
$\ZZ_p$-module $\Lambda$, we have a canonical isomorphism
$M^\sigma \cong M$ which we use to identify these two $\What(R)$-modules.
\end{segment}

\begin{segment}\label{defsigma}
Let $X$ be any $\What(R)$-scheme. Then the ring endomorphism $\sigma$ of
$\What(R)$ induces a map $\sigma\colon X(\What(R)) \to
X(\What(R))$. We will use this notation in particular for the group
scheme $X = \Gscr \otimes_{\ZZ_p} \What(R)$ and for the scheme of parabolics of
$\Gscr\otimes_{\ZZ_p} \What(R)$.
\end{segment}

\begin{segment}\label{defdisplay}
\begin{segproc}{Definition} We set $M :=
  \Lambda \otimes_{\ZZ_p} \What(R)$. Then $M$ carries a symplectic
  form $\lrangle$. A \emph{split
  symplectic Dieudonn\'e display over $R$} consists of a tuple
  $(S,T,F,V^{-1})$ where $S$ and $T$ are totally isotropic
  $\What(R)$-submodules of $M$ such that $S \oplus T = M$. Further
  $F\colon M^{\sigma} \to M$ and $V^{-1}\colon Q^\sigma \to M$ are
  $\What(R)$-linear maps where $Q := S \oplus \Ihat_RT = S
  + \Ihat(R)M \subset M$. The following properties are satisfied:
\begin{definitionlist}
\item
$V^{-1}$ is surjective.
\item
For all $x \in M$ and $w \in \What(R)$ we have
\[
V^{-1}(1 \otimes \tau(w)x) = wF(1 \otimes x).
\]
\item
For all $y, y'\in Q$ we have
\[
\tau(\langle V^{-1}(1\otimes y), V^{-1}(1 \otimes y') \rangle) =
\langle y, y' \rangle.
\]
\end{definitionlist}
\end{segproc}
\end{segment}

\begin{segment}\label{remarkdisplay}
By conditions (a) and (b), the map
\[
V^{-1} \oplus F\colon S^\sigma \oplus T^\sigma \rTo M
\]
is a surjective $\What(R)$-linear map of free
$\What(R)$-modules of the same rank, hence it is an isomorphism.

As $S$ and $T$ are both totally isotropic and as $M = S \oplus T$, we
have $\rk_{\What(R)}(S) = \rk_{\What(R)}(T) = g$.
\end{segment}

\begin{segment}\label{remarkduality}
Using the identity $F(1 \otimes x) = V^{-1}(1 \otimes \tau(1)x)$ it is
easy to check that \eqref{defdisplay} implies
\begin{alignat*}{2}
\langle V^{-1}(1 \otimes y), F(1 \otimes x)\rangle &= \sigma(\langle
y,x \rangle) & & \quad\text{for all $y \in Q$ and $x \in P$,}\\
\langle F(1 \otimes x), F(1 \otimes x') \rangle &= p\sigma(\langle
x,x' \rangle) & & \quad\text{for all $x,x' \in P$.}
\end{alignat*}
\end{segment}

\begin{segment}\label{associatedisplay}
Let $(A,\lambda,\alpha)$ be an
$R$-valued point of $\Ascr^\#_0$ \eqref{covermoduli}. Let
$(M,Q,F,V^{-1})$ be the Dieudonn\'e display associated to the
$p$-divisible group of $A$ by the theory of Zink \cite{Zi}.

Moreover, $\lambda$ induces a perfect alternating form $\lrangle$ on
the free $\What(R)$-module $M$ such that $\tau\langle V^{-1}(1
\otimes y), V^{-1}(1 \otimes y') \rangle = \langle y, y' \rangle$
for all $y, y' \in Q$.

By \eqref{localisom} we can find a $\What(R)$-linear
symplectic isomorphism
\[
\agtilde\colon M \rTo^{\sim} \Lambda \otimes \What(R)
\]
whose reduction modulo $\Ihat(R)$ is equal to $\alpha$.
We use $\agtilde$ to identify $M$ and $\Lambda \otimes \What(R)$.

By definition of a Dieudonn\'e display we have a split exact
sequence of free $R$-modules
\[
0 \rTo C \rTo M/\Ihat_R M \rTo M/Q \rTo 0.
\]
We choose a totally isotropic direct summand $S$ of $M$
such that its reduction modulo $\Ihat_R$ is equal to $C$ and we choose a
totally isotropic $\What(R)$-complement $T$ of $S$ in
$M$. Then $(S,T,F,V^{-1})$ is a split symplectic Dieudonn\'e display such
that $Q = S + \Ihat(R)M$.
\end{segment}

\begin{segment}\label{groupdisplay}
We are now going to give a group theoretic reformulation of the
split Dieudonn\'e displays with additional structures defined in
\eqref{defdisplay}: For any ring $R$ as in \eqref{defWhat} we define
$\Ytilde_J(\What(R))$ to be the set of triples
$(\Ptilde,\Mtilde,\gtilde)$ where $\Ptilde \subset \Gscr'_{\What(R)}$
is a parabolic of type $J$, where $\Mtilde \subset \Ptilde$ is a Levi
subgroup and where $\gtilde \in \Gscr'(\What(R))$.
\end{segment}

\begin{segment}\label{displaytogroup}
Let $(S,T,F,V^{-1})$ be a split symplectic Dieudonn\'e display over $R$.
We associate an element $(\Ptilde,\Ltilde,\gtilde)$ in
$\Ytilde_J(\What(R))$ as follows: We define $\Ptilde$ as the
stabilizer of the flag $0 \subset S \subset M$ in
$\Gscr'_{\What(R)}$. Then $\Ptilde$ is a parabolic of type $J$ by
\eqref{defdisplay}~(d).

Furthermore, $\Ltilde$ is by definition the stabilizer of the
decomposition $M = S \oplus T$ in $\Gscr'_{\What(R)}$. Clearly, $\Ltilde$ is a
Levi subgroup of $\Ptilde$.

Finally let $\gtilde$ be the composition
\[
M \rTo^{\sim} M^\sigma = S^\sigma \oplus T^\sigma \rTo^{V^{-1}|S^\sigma
  \oplus F|T^\sigma} M.
\]

\begin{proclamation}{Lemma}
The map constructed above defines a bijection between the set of all
split symplectic Dieudonn\'e displays over $R$ and the set
$\Ytilde_J(\What(R))$.
\end{proclamation}

\begin{proof}
Clearly $\gtilde$ is a $\What(R)$-linear map. By
\eqref{remarkdisplay} it is an isomorphism. Now we use
\eqref{remarkduality} to check that $\gtilde$ respects the
alternating form $\lrangle$ (and therefore $\gtilde \in
\Gscr'(\What(R))$):

Let $t,t' \in T$ and write $t = \sum_i w_i \otimes \lambda_i$ and $t'
= \sum_j w'_j \otimes \lambda'_j$ with $w_i,w'_j \in \What(R)$ and
$\lambda_i, \lambda'_j \in \Lambda$. As $T$ is totally isotropic, we
have
\begin{align*}
\langle t,t' \rangle &= 0.\\
\intertext{On the other hand}
\langle \gtilde(t), \gtilde(t') \rangle &= \sum_{i,j}w_iw'_j\langle
F(1 \otimes (1 \otimes \lambda_i)), F(1 \otimes (1 \otimes
\lambda'_j)) \rangle\\
&= \sum_{i,j}w_iw'_jp\sigma(\langle 1 \otimes \lambda_i, 1 \otimes
\lambda'_j\rangle)\\
&= p\langle \sum_i w_i \otimes \lambda_i, \sum_j w'_j \otimes
\lambda'_j \rangle\\
&= 0.
\end{align*}
A similar argument shows that $\langle \gtilde(s), \gtilde(s') \rangle
= \langle s,s' \rangle$ for all $s,s' \in S$. For $s = \sum_i w_i
\otimes \lambda_i \in S$ and $t = \sum_j w'_j \otimes \lambda'_j \in
T$ we have
\begin{align*}
\langle \gtilde(s), \gtilde(t) \rangle &= \sum_{i,j}w_iw'_j\langle
V^{-1}(1 \otimes (1 \otimes \lambda_i)), F(1 \otimes (1 \otimes
\lambda'_j)) \rangle\\
&= \sum_{i,j}w_iw'_j\sigma(\langle 1 \otimes \lambda_i, 1 \otimes
\lambda'_j\rangle)\\
&= \langle \sum_i w_i \otimes \lambda_i, \sum_j w'_j \otimes
\lambda'_j \rangle\\
&= \langle s,t \rangle.
\end{align*}
This shows that $(\Ptilde,\Ltilde,\gtilde) \in \Ytilde_J(\What(R))$.

We construct an inverse map: Let $(\Ptilde,\Ltilde,\gtilde)$ be in
$\Ytilde_J(\What(R))$. We let $S$ be the unique direct summand of $M$
such that its stabilizer is equal to $\Ptilde$ and let $T \subset
M$ be the unique direct complement of $S$ such that the
stabilizer of the decomposition $S \oplus T$ is equal to
$\Ltilde$. Further we set for $t \in T$, $s \in S$, $w \in \What(R)$
\begin{align*}
F(1 \otimes t) &= \gtilde(t), & F(1 \otimes s) &= p\gtilde(s), \\
V^{-1}(1 \otimes \tau(w)t) &= w\gtilde(t), & V^{-1}(1 \otimes s) &= \gtilde(s),
\end{align*}

Clearly $V^{-1}$ is surjective as $\gtilde$ is surjective. Moreover,
we have for $w \in \What(R)$, $t \in T$ and $s \in S$
\begin{gather*}
V^{-1}(1 \otimes \tau(w)t) = w\gtilde(t) = wF(t),\\
V^{-1}(1 \otimes \tau(w)s) = V^{-1}(pw \otimes s) = pw\gtilde(s) =
wF(1 \otimes s)
\end{gather*}
which shows that condition (b) of \eqref{defdisplay} holds.
A similar although much more lengthy calculation shows that
condition~(c) is also satisfied.

This shows that $(S,T,F,V^{-1})$ is a split
symplectic Dieudonn\'e display. Clearly this construction defines an inverse.
\end{proof}
\end{segment}


\section{Flatness of the mod $p$ period morphism}

\begin{segment}
Let $w_0$ be the element of maximal length in $W$ and let $x$ be the element
of minimal length in $W_Jw_0W_J = w_0W_J$.

We denote by $Z_J$ the functor on $\FF_p$-schemes which is the
Zariski-sheafi\-fication of the functor $Z^{\prime}_J$ which associates to an
$\FF_p$-scheme $S$ the set of triples
\[
(P,Q,U_QgU_{F(P)})
\]
where $P \subset
\Gbar_S$ and $Q \subset \Gbar_S$ are parabolics of type $J$, and where
$g \in G(S)$ is an element such that $\relpos(Q,{}^gF(P)) = x$. By
\cite{MW}~3.12 this functor is representable by a scheme. For any
affine scheme $S$ we have $Z_J(S) = Z^{\prime}_J(S)$.

By definition of $x$ we have $\relpos(Q,{}^gF(P)) = x$ if and only
if\break $Q \cap {}^gF(P)$ is a common Levi subgroup of $Q$ and
${}^gF(P)$, i.e.~$Q$ and ${}^gF(P)$ are in opposition.
\end{segment}

\begin{segment}\label{torsoroverPar}
The forgetful morphism $Z_J \to \Par_J \times \Par_J$ which is defined
on points by $(P,Q,U_QgU_{F(P)}) \mapsto (P,Q)$ makes $Z_J$ into a
torsor over $\Par_J \times \Par_J$ under a reductive group scheme of
dimension $\dim(P/U_P)$ for any parabolic $P$ of $G$ of type $J$
(\cite{MW}~3.11). In particular, $Z_J$ is a smooth $\FF_p$-scheme whose
dimension equals $\dim(G)$.
\end{segment}

\begin{segment}
We denote by $\Ztilde_J$ the scheme which represents the functor on
$\FF_p$-schemes which associates to $S$ the set of triples $(P,Q,g)$ where
$P \subset \Gbar_S$ and $Q \subset \Gbar_S$ are parabolics of type $J$
and where $g \in G(S)$ is an element such that $\relpos(Q,{}^gF(P)) = x$.

Via the forgetful morphismus $(P,Q,g) \mapsto (P,Q)$ we will consider
$\Ztilde_J$ as a scheme over $\Par_J \times \Par_J$.

We consider $\Uscr_J \times F(\Uscr_J)$ as a group scheme over $\Par_J \times
\Par_J$. Then $(u,v)\cdot(P,Q,g) = (P,Q,ugv^{-1})$ defines an action
of $\Uscr_J \times F(\Uscr_J)$ on $\Ztilde_J$ over $\Par_J \times
\Par_J$. The fppf-quotient of this action is $Z_J$.
\end{segment}

\begin{segment}\label{Ztildetorsor}
\begin{segproc}{Lemma} The action of $\Uscr_J \times F(\Uscr_J)$ on
  $\Ztilde_J$ is free and hence $\Ztilde_J$ is a torsor under the
  smooth group scheme $\Uscr_J \times F(\Uscr_J)$ over $Z_J$.
\end{segproc}
\begin{proof}
Let $(P,Q,g) \in \Ztilde_J$ and let $u \in U_Q$ and $v \in U_{F(P)}$
such that $ugv^{-1} = g$. This implies that $gv^{-1}g^{-1} = u \in U_Q
\cap {}^gU_{F(P)}$. But by definition of $\Ztilde_J$, ${}^gF(P)$ and
$Q$ are in opposite position. Therefore, $U_Q \cap {}^gU_{F(P)} = (1)$
which implies $u = v = 1$.
\end{proof}
\end{segment}

\begin{segment}
The group $\Gbar$ acts on $\Ztilde_J$ by the rule
$$
h\cdot(P,Q,g) = ({}^hP,{}^hQ,hgF(h)^{-1}).
$$
This induces an action of $\Gbar$ on $Z_J$.
\end{segment}

\begin{segment}
We denote by $\Ytilde_J$ the $\FF_p$-scheme which represents the functor
which associates to every $\FF_p$-scheme $S$ the set of triples $(P,L,g)$
where $P \in \Par_J(S)$, $g \in \Gbar(S)$ and $L$ is a Levi subgroup of
$P$.

We construct a morphism
\[
\Ytilde_J \rTo \Ztilde_J.
\]
For every $S$-valued point $(P,L,g)$ of $\Ytilde_J$ we
 define $Q$ as the unique parabolic of type $J$ such that
 ${}^{g^{-1}}Q \cap F(P) = F(L)$ (\cite{SGA3}~Exp.~XXVI,~4.3.). Then
 $(P,Q,g) \in \Ztilde_J(S)$.
\end{segment}

\begin{segment}
Now we relate the moduli spaces defined in \eqref{covermoduli} and the
varieties $\Ytilde_J$ and $Z_J$.

We define a morphism
$$
\zgtilde\colon \Ascrtilde_0 \to \Ytilde_J
$$
as follows: To every $S$-valued
point $(A,\lambda,\alpha,C',D')$ we associate the triple
$(P,L,g)$ where $P$ is the stabilizer of $\alpha(C)$ in $\Gbar_S$,
where $L$ is the stabilizer of the decomposition $\alpha(C) \oplus
\alpha(C') = \Lambda_S$, and where $g$ is the
composition
\[
\Lambda_S \rTo^{\sim} \Lambda_S^{(p)} = \alpha(C)^{(p)} \oplus
\alpha(C')^{(p)} \rTo^{\sim}_{\varphi_1 \oplus \varphi_0} \alpha(D)
\oplus \alpha(D') = \Lambda_S.
\]

By definition $L$ is a Levi subgroup of $P$, hence $(P,L,g) \in
\Ytilde_J(S)$.
\end{segment}

\begin{segment}\label{open}
It follows from the definitions that $\zgtilde$ induces a morphism
$$
\zeta\colon \Ascr^{\#}_0 \to Z_J.
$$

\begin{proclamation}{Theorem}
The morphism $\zeta$ is flat.
\end{proclamation}
\end{segment}

\begin{segment}\label{opencrit}
We will show that $\zeta$ is universally open. As $Z_J$ and
$\Ascr^{\#}_0$ are both regular, it then follows from
\cite{EGA}~IV,~(15.4.2), that $\zeta$ is flat. To show that $\zeta$ is
universally open, we use the following criterion.

\begin{proclamation}{Proposition}
Let $Y$ be a noetherian geometrically unibranch scheme, let $X$ be a
scheme, and let $f\colon X \to Y$ be a morphism of finite
type. Assume that for every commutative diagram
\begin{equation}\label{valdiagram}
\begin{diagram}
\Spec(k) & \rTo & X \\
\dTo & & \dTo^{f} \\
\Spec(R) & \rTo^g & Y \\
\end{diagram}
\end{equation}
where $R$ is a complete discrete valuation ring with algebraically closed
residue field $k$, there
exists a surjective morphism $\Spec(\Rtilde) \to \Spec(R)$ of discrete
valuation rings and a morphism $\gtilde\colon \Spec(\Rtilde) \to X$ which
commutes with
\[
\begin{diagram}
\Spec(\kappa(\Rtilde)) & \rTo & \Spec(k) & \rTo & X \\
\dTo & & & & \dTo^{f} \\
\Spec(\Rtilde) & \rTo & \Spec(R) & \rTo^{g} & Y. \\
\end{diagram}
\]
Then $f$ is universally open.
\end{proclamation}

\begin{proof}
By \cite{EGA} IV, 14.4.1 it suffices to show that $f$ is open.
Let $U \subset X$ be an open subset. By Chevalley's theorem $f(U)$ is
contructible. By \cite{AK} V, 4.4 it suffices therefore to show that
$f(U)$ is stable under generization. Let $x_0 \in U$ and $y_0 = f(x_0)
\in f(U)$ and let $y_1 \in Y$ be a generization with $y_1 \not= y_0$.

By \eqref{existval} below, there exists a diagram like in
\eqref{valdiagram}. We apply the hypothesis and find a morphism
$\gtilde\colon \Spec(\Rtilde) \to X$ such that $f \circ \gtilde$ is
the composition
\[
\Spec(\Rtilde) \rTo \Spec(R) \rTo \Spec(R') \rTo Y.
\]
The image $x_1$ of the generic point of $\Spec(\Rtilde)$ under $\gtilde$ is a
generization of $x_0$ and hence lies in $U$ as $U$ is open, and therefore $y_1
= f(x_1) \in f(U)$. 
\end{proof}
\end{segment}

\begin{segment}\label{existval}
\begin{segproc}{Lemma}
Let $Y$ be a locally noetherian scheme and let $f\colon X \to Y$ be a
morphism of schemes. Let $x_0 \in X$, $y_0 := f(x_0)$ and let $y_1
\not= y_0$ be a generization of $y_0$. Then there exists a commutative diagram
\[
\begin{diagram}
\Spec(\kappa) & \rTo^{h} & X \\
\dTo_{i} & & \dTo^{f} \\
\Spec(R) & \rTo^{g} & Y \\
\end{diagram}
\]
where $R$ is a discrete valuation ring with algebraically closed residue field
$i\colon \Spec(\kappa) \air \Spec(R)$ such that the image of the generic
(resp.\ special) point of $\Spec(R)$ under $g$ is $y_1$ (resp.\ $y_0$) and
such that the image of $h$ is $x_0$.
\end{segproc}

\begin{proof}
There exists a morphism
$g'\colon \Spec(R') \to Y$ where $R'$ is a discrete valuation ring such that
$g'(s') = y_0$ and $g'(\eta') = y_1$ where $s'$ (resp.\ $\eta'$) is the closed
(resp.\ generic) point of $\Spec(R')$.

Let $\mfr'$ be the maximal ideal of $R'$ and let $\kappa$ be
an algebraically closed field extension of $\kappa(y_0)$ such that there exist
$\kappa(y_0)$-embeddings $\kappa(x_0) \air \kappa$ and $\kappa(s') \air
\kappa$ and let $R' \to R$ be a flat local homomorphism of $R'$ into a
complete discrete valuation ring $R$ with residue field $\kappa$ such that
$\mfr'R$ is the maximal ideal of $R$ (this exists by \cite{EGA} ${\bf
  0}_I$ 6.8.3). We set $g$ as the composition
\[
\Spec(R) \rTo \Spec(R') \rTo^{g'} Y
\]
and $h$ as the composition
\[
\Spec(\kappa) \rTo \Spec(\kappa(x_0)) \rTo X.
\]
\end{proof}
\end{segment}

\begin{segment}\label{proofofopen} {\it Proof of \eqref{open}}.
By \eqref{opencrit} it suffices to show the following lemma.

\begin{proclamation}{Lemma}
Let $k$ be an algebraically closed field of characteristic $p$, let\break $R
= k\dlbrack \eps \drbrack$ be the ring of formal power series in one
variable $\eps$ and set $R_1 = k\dlbrack \eps^{1/p} \drbrack$. We denote by
\[
h\colon \Spec(R_1) \rTo \Spec(R)
\]
the natural morphism.

Let $x = (A,\lambda,\alpha)$ be a $k$-valued point of
$\Ascr^\#_0$. Let $(P,Q,[g]) \in Z_J(k)$ be the image of $x$
under $\zeta$. Denote by
$(P_\eps,Q_\eps,[g_\eps]) \in Z_J(R)$ any deformation of $(P,Q,g)$ to
$R$. Then there exists a deformation $x_1 =
(A_1,\lambda_1,\iota_1,\alpha_1) \in \Ascr^\#_0(R_1)$ of $x$ such that
$h(\zeta(x_1)) = (P_\eps,Q_\eps,[g_\eps])$.
\end{proclamation}

\begin{proof}
Let $\Pscr = (M,Q,F,V^{-1})$ be the Dieudonn\'e display of the $p$-divisible
group of the abelian variety $A$. The free $W(k)$-module $M$ is
equipped with a perfect alternating form via $\lambda$. Moreover, we can fix an
identification $\agtilde$ of $M$ as a symplectic $W(k)$-module with
$\Lambda \otimes W(k)$ which lifts the isomorphism $\alpha$ \eqref{localisom}.
Set $M_\eps = \Lambda \otimes \What(R)$ and $M_{\eps_1} = \Lambda
\otimes \What(R_1)$.

We choose submodules $S \subset M$ and $T \subset
M$ such that $(S,T,F,V^{-1})$ is a split symplectic Dieudonn\'e
display over $k$. Let $(\Ptilde,\Ltilde,\gtilde)$ be the associated
element in $\Ytilde_J(W(k))$ \eqref{displaytogroup}.

Let $g \in [g]$ be the reduction of $\gtilde$ modulo $p$ and choose
$g_\eps \in [g_\eps]$ such that the reduction of $g_\eps$ modulo
$\eps$ is equal to $g$. Set $P_{\eps,1} := P_\eps \otimes_R R_1$,
$Q_{\eps,1} := P_\eps \otimes_R R_1$ and let $g_{\eps,1}$ be the
element $g_\eps$ considered as an $R_1$-valued point of $\Gscr'$.

Let $L_{\eps,1}$ be a Levi subgroup of
$P_{\eps,1}$ such that
\[
F(L_{\eps,1}) = {}^{(g^{-1}_{\eps,1})}Q_{\eps,1} \cap F(P_{\eps,1}).
\]

For any smooth $W(k)$-scheme $X$ the canonical map
\[
X(\What(R_1)) \to X(R_1) \times_{X(k)} X(W(k))
\]
is surjective. Applying this to the
scheme of parabolic subgroups of $\Gscr_{W(k)}$ of type $J$, to the group
scheme $\Gscr_{W(k)}$ itself and to the smooth schemes of Levi subgroups of
fixed parabolic subgroup of $\Gscr_{W(k)}$, we see that
there exists an element
\[
(\Ptilde_{\eps,1},\Ltilde_{\eps,1},\gtilde_{\eps,1}) \in
 \Ytilde_J(\What(R_1))
\]
whose reduction to $R_1$ equals $(P_{\eps_1},L_{\eps_1},g_{\eps,1})$
and whose reduction to $W(k)$ equals $(\Ptilde,\Ltilde,\gtilde)$.

Let $(S_{\eps,1},T_{\eps,1},F_{\eps,1},V^{-1}_{\eps,1})$ be the split
symplectic Dieudonn\'e display associated to
$(\Ptilde_{\eps,1},\Ltilde_{\eps,1},\gtilde_{\eps,1})$
\eqref{displaytogroup}. We set $Q_{\eps,1} = S_{\eps,1} \oplus \Ihat_R
T_{\eps,1}$. Then $(M_{\eps,1}, Q_{\eps,1},
F_{\eps,1},V^{-1}_{\eps,1})$ is a Dieudonn\'e display which lifts
$\Pscr$. Via the correspondence of $p$-divisible groups over $R_1$ and
Diedonn\'e displays over $R_1$ \cite{Zi}, the theorem follows from
Serre-Tate theory.
\end{proof}
\end{segment}

\begin{segment}
The morphism $\zeta\colon \Ascr^\#_0 \rTo Z_J$ induces a flat
morphism
\[
\zgbar\colon \Ascr_0 \to [G \backslash Z_J]
\]
where $[G\backslash Z_J]$ denotes the stack quotient of $Z_J$ by the action
of $G$. For every $G$-orbit $O$ of $Z_J$, the quotient stack $[G
\backslash O]$ is a locally closed substack of $[G\backslash Z_J]$
whose underlying topological space consists of only one point. From
\cite{MW}~3.25 we know that the $G$-orbits of $Z_J$ are in natural
bijection to ${}^JW$. Moreover, we know that the $G$-orbit $O^u$
corresponding to $u \in {}^JW$ has codimension $\dim(\Par_J) -
\ell(u)$ in $Z_J$. Therefore, the same holds for the substack $[G
\backslash O^u]$ of $[G \backslash Z_J]$.

The inverse image of $[G \backslash O^u]$ in $\Ascr_0$ is denoted by
$\Ascr_0^u$. These are just the Ekedahl-Oort strata.
\end{segment}

\begin{segment}\label{dimension}
\begin{segproc}{Corollary}
The $\Ascr_0^u$ for $u \in {}^JW$ form a stratification of $\Ascr_0$
(i.e.~they are locally closed and the closure of one stratum is the
union of strata). All strata $\Ascr^u_0$ are equi-dimensional and we have
\[
\dim(\Ascr_0^u) = \ell(u).
\]
\end{segproc}

\begin{proof}
Clearly the $O^u$ (for $u \in {}^JW$) form a stratification of $Z_J$
as they are just the $G$-orbits. As $\zgbar$ is open, this is true for the
$\Ascr^u_0$ as well.

Moreover, the unique closed point of the underlying topological space
of $[G'\backslash Z_J]$ is contained in the image of $\zgbar$ (any
superspecial principally polarized abelian variety is mapped to this
point). Therefore the openness of $\zgbar$ implies that $\zgbar$ is
surjective. In other words, all strata $\Ascr^u_0$ are nonempty.

As $\zgbar$ is flat, it also respects
codimension, and hence we have
\begin{align*}
\dim(\Ascr_0^u) &= \dim(\Ascr_0) - \codim(\Ascr_0^u,\Ascr_0) \\
&= \dim(\Par_J) - (\dim(\Par_J) - \ell(u)) \\
&= \ell(u).
\end{align*}
\end{proof}
\end{segment}

\begin{segment}
\begin{segproc}{Corollary}
The morphism
\[
\zgbar\colon \Ascr_0 \to [G \backslash Z_J]
\]
is faithfully flat.
\end{segproc}
\end{segment}

\begin{segment}
The dimension formula \eqref{dimension} has also been shown by Oort in
\cite{Oo}. Here we give a new proof which can be carried over to
arbitrary good reductions of Shimura varieties of PEL-type. We will
come back to this in \cite{We3}.
\end{segment}

\begin{segment}
We can identify ${}^JW$ as a set with $\{0, 1\}^g$
\eqref{Siegelparabolic}. Further we know by \cite{We2} that all
Ekedahl-Oort strata are nonempty (although this was certainly known
before). Now the corollary \eqref{dimension} tells us that the
Ekedahl-Oort stratum corresponding to $u =
\{\epsilon_1,\dots,\epsilon_g\} \in {}^JW$ is equidimensional of
dimension $g\epsilon_1 + (g-1)\epsilon_2 + \dots + \epsilon_g$.
\end{segment}

\begin{segment}\label{modelsmooth}
Now consider the forgetful morphism
$\eta\colon Z_J \to \Par_J$ which is defined on points by $(P,Q,g)
\mapsto P$. We know from \eqref{torsoroverPar} that $\eta$ is smooth.
Moreover, by Grothendieck-Messing theory we know that the composition
$\theta := \eta \circ \zeta\colon \Ascr^\#_0 \rTo \Par_J$ is
smooth. We therefore have a diagram of morphisms
\[
\begin{diagram}
\Ascr^\#_0 & & \rTo^{\zeta} & & Z_J \\
& \rdTo_\theta & & \ldTo_{\eta} & \\
& & \Par_J, & & \\
\end{diagram}
\]
where $\eta$ and $\theta$ are smooth and where $\zeta$ is flat.
\end{segment}



\begin{thebibliography}{SGA3}
\bibitem[AK]{AK} A.~Altman, S.~Kleiman: {\it Introduction to Grothendieck
Duality Theory}, Lecture Notes in Mathematics {\bf 146}, Springer, 1970.
\bibitem[EGA]{EGA} A.~Grothendieck, J.~Dieudonn\'e: {\it El\'ements de
    g\'eom\'etrie alg\'ebrique}, Publ.\ Math.\ de l'IHES {\bf 4}, {\bf 8}, 
{\bf 11}, {\bf 17}, {\bf 20}, {\bf 24}, {\bf 28}, {\bf 32} 
(1960--67).
\bibitem[Lu]{Lu} G.~Lusztig: {\it Cuspidal Local systems and Graded
    Hecke Algebras II}, in \cite{RoG}, pp. 217--275.
\bibitem[MAV]{MAV} C.~Faber, G.\ van der Geer, F.~Oort (eds.): {\it 
Moduli of abelian varieties\/}, Progr.\ in Math.\ {\bf 195}, 
Birkh\"auser Verlag, Basel, 2001.
\bibitem[Mo1]{Mo1} B.~Moonen: {\it Group schemes with additional 
structures and
Weyl group cosets\/}, in: \cite{MAV}, pp.~255--298.
\bibitem[Mo2]{Mo2} B.~Moonen: {\it A dimension formula for Ekedahl-Oort
strata\/}, preprint 2002, math.AG/0208161, to appear in
the Ann.\ Inst.\ Fourier (Grenoble).
\bibitem[MW]{MW} B.~Moonen, T.~Wedhorn: {\it Discrete invariants of
    varities in positive characteristic}, preprint Bonn 2004, to
    appear in Int. Math. Res. Not.
\bibitem[Oo]{Oo} F.~Oort: {\it A stratification of a moduli space of 
abelian varieties\/}, in: \cite{MAV}, pp.~435--416.
\bibitem[RoG]{RoG} B.~Allison, G.~Cliff (eds.): {\it Representations
    of Groups}, CMS Conference Proceedings {\bf 16}, Canadian Mathematical
    Society, 1995.
\bibitem[RZ]{RZ} M.~Rapoport, T.~Zink: {\it Period spaces for 
$p$-divisible groups\/},
Annals of Math.\ Studies {\bf 141}, 
Princeton Univ.\ Press, Princeton, NJ, 1996.
\bibitem[SGA3]{SGA3} M.~Demazure et al.: {\it Sch\'emas en groupes, I,
II, III\/}, LNM {\bf 151}, {\bf 152}, {\bf 153}, 
Springer-Verlag, Berlin, 1970.
\bibitem[Wd1]{We1} T.~Wedhorn: {\it Ordinariness in good reductions of
    Shimura varieties of PEL-type}, Ann.\ scient.\ \'Ec.\ Norm.\ 
Sup.\ (4), {\bf 32}, (1999), 575--618.
\bibitem[Wd2]{We2} T.~Wedhorn: {\it The dimension of Oort strata of 
Shimura varieties of PEL-type}, in: \cite{MAV}, pp.~441--471.
\bibitem[Wd3]{We3} T.~Wedhorn: {\it $F$-Zips strata for Shimura
    varieties of PEL type}, in preparation
\bibitem[Zi]{Zi} T.~Zink: {\it A Dieudonn\'e Theory for $p$-divisible
    groups},in: Class Field Theory, Its Centenary and Prospect
    pp.~1--22, Advanced Studies in Pure Mathematics, Tokyo 2000
\end{thebibliography}
\end{document}